\documentclass[11pt]{article}
\usepackage{amsmath} 
\usepackage{amsthm}
\usepackage{geometry}                
\usepackage[parfill]{parskip}    
\usepackage{graphicx}
\usepackage[noend]{algorithm2e}
\usepackage{amssymb}
\usepackage{epstopdf}
\usepackage{verbatim}
\usepackage{hyperref}
\hypersetup{colorlinks,%
citecolor=black,%
filecolor=black,%
linkcolor=black,%
urlcolor=black%
}
\newtheorem{thm}{Theorem}
\newtheorem*{thma}{Theorem 1}

\newtheorem{cor}[thm]{Corollary}
\newtheorem{lem}[thm]{Lemma}

\theoremstyle{definition}
\newtheorem{defn}[thm]{Definition}

\theoremstyle{remark}
\newtheorem{rem}[thm]{Remark}
\newcommand{\f}[1]{\mathbb{F}_{#1}}
\newcommand{\tfr}[1]{{\widehat{#1}}}
\newcommand{\K}{\mathcal{K}}
\newcommand{\ZZ}{\mathbb{Z}}
\newcommand{\QQ}{\mathbb{Q}}

\DeclareMathOperator{\Tr}{Tr}
\DeclareMathOperator{\wt}{wt}
\DeclareMathOperator{\gal}{Gal}

\begin{document}
\title{Some Congruences of Kloosterman Sums and their Minimal Polynomials}
\author{Faruk G\"olo\u glu\thanks{Research supported by Claude Shannon Institute,
Science Foundation Ireland Grant 06/MI/006}, Gary McGuire\footnotemark[1]
, Richard Moloney\thanks{Research supported by Claude Shannon Institute,
Science Foundation Ireland Grant 06/MI/006, 
and the Irish Research Council
for Science, Engineering and Technology}\\
School of Mathematical Sciences\\
University College Dublin\\
Ireland}

\date{}

\maketitle

\begin{abstract}
We prove two results on Kloosterman sums over finite fields, using Stickelberger's theorem
and the Gross-Koblitz formula.  The first result concerns the minimal polynomial 
over $\mathbb{Q}$ of a
Kloosterman sum, and the second result gives
a characterisation of ternary Kloosterman sums modulo 27.
\end{abstract}

\section{Introduction}\label{Intro}
Let $p$ be an odd prime, $n \ge 1$ an integer, $q = p^n$ and $\zeta$ a primitive $p^{\text{th}}$ root of unity. We let $\f{q}$ denote the finite field with $q$ elements, and let 
$\Tr$ denote the absolute trace function $\Tr:\f{q} \to \f{p}$, \[
\Tr (a) = a + a^p + a^{p^2} + \dots + a^{p^{n-1}}\,.
\] The Kloosterman sum of $a \in \f{q}$ is defined to be \[\K_q(a) = \sum_{x \in \f{q}}\zeta^{\Tr(x^{-1}+ax)}
\]
where we interpret $0^{-1}$ as $0$.
We remark that some authors do not include $0$ in the definition of Kloosterman sum.


Obviously $\K_{q} (a)$ is an algebraic integer lying in the 
cyclotomic field $\mathbb{Q}(\zeta)$. 
It is well known  that
\[\gal (\mathbb{Q}(\zeta) / \mathbb{Q}) = \{ \zeta \mapsto \zeta^{i} \ |\ i \in (\mathbb{Z}/p\ZZ)^{*}\},\] and 
it is easy to show (see \cite{finns}) that 
the Galois automorphism $ \zeta \mapsto \zeta^{i}$ has the effect
$\K_q(a) \mapsto \K_q(i^2 a)$, for any integer $i$.
If we let
\[c_a(x) = \prod_{i=1}^{\frac{p-1}{2}}(x-\K_q(i^2a))\] 
it follows that $c_a(x)$ (which has degree $(p-1)/2$) is
the characteristic polynomial of $\K_q(a)$ over $\QQ$. 
If $m_a(x)$ is the minimal polynomial of $\K_q(a)$ over $\QQ$, then
\[c_a(x) = m_a(x)^{e_a}\] for some $e_a$ dividing $\frac{p-1}{2}$.
Most of the time, it is true that $e_a=1$.
For example, Wan \cite{Wan} showed that $e_a=1$ if $\Tr(a)\not=0$.

Moisio \cite{CertainValues} considered the reduction of the minimal polynomial $m_a(x)$ modulo $p$.
He showed that all coefficients, apart from the leading coefficient, are divisible by $p$.

In this paper, our first result concerns   the reduction of the minimal polynomial $m_a(x)$ modulo $p^2$.
 In Section \ref{MainProof}, we prove the following result about the constant term. 
\bigskip

\begin{thm}\label{main}
Let $p$ be an odd prime, and let $\left(\frac{\cdot}{p}\right)$ be the Legendre symbol. Then 
\[\prod_{i=1}^{\frac{p-1}{2}}\K_q(i^2a) \equiv p \left(\frac{\Tr(a)}{p}\right) \pmod{p^2}\,.\]
\end{thm}

As a consequence, the constant term of the characteristic polynomial, 
which is \[(-1)^{\frac{p-1}{2}}\prod_{i=1}^{\frac{p-1}{2}}(\K_q(i^2a)),
\]
is always congruent to either 0 or $\pm p$ mod $p^2$.

In the case that $p = 3$, Theorem \ref{main} becomes the following theorem.

\begin{thm}\label{mainternary}
Let $n>1$. For $a \in \f{3^n}$,
\[
\mathcal{K}_{3^n}(a) \equiv \left\{ \begin{array}{ll} 
0 \ (\mathrm{mod} \ 9) & \textrm{if } \Tr(a) = 0,\\
3 \ (\mathrm{mod} \ 9) & \textrm{if } \Tr(a) = 1,\\ 
6 \ (\mathrm{mod} \ 9) & \textrm{if } \Tr(a) = 2. 
\end{array}
\right. 
\]
\end{thm}

This is precisely the modulo 9 characterisation of the ternary Kloosterman sum which 
we previously proved in \cite{SetaFGR}.
The second result of this paper, 
see Corollary \ref{main27} in Section \ref{GKmod27}, is to 
extend this result to a modulo 27 characterisation of the ternary Kloosterman sum. 

\section{Background}

In this section we present the background information that is used in our proofs.

\subsection{Teichm\"uller characters and Gauss sums}

Consider multiplicative characters taking their values in an algebraic extension of $\mathbb
{Q}_p$. Let $\xi$ be a primitive $(q-1)^{\text{th}}$ root of unity in a 
fixed algebraic closure of  $\mathbb{Q}_p$.  
The group of multiplicative characters of $\mathbb{F}_q$ (denoted $
\widehat{\mathbb{F}_q^{\times}}$) is cyclic of order $q-1$. 
The group $\widehat{\mathbb{F}_q^
{\times}}$ is generated by the Teichm\"uller character $\omega : \mathbb{F}_q \to 
\mathbb{Q}_p(\xi)$, which, for a fixed generator $t$ of $\mathbb{F}_q^{\times}$, is defined 
by $\omega (t^j) = \xi^j$. We set $\omega (0)$ to be $0$. An equivalent definition is that $\omega$ satisfies \[\omega(a) \equiv a \pmod{p}\] for all $a \in \f{q}$.

Let $\zeta$ be a fixed primitive $p$-th root of unity in the 
fixed algebraic closure of  $\mathbb{Q}_p$.  
Let $\mu$ be the canonical additive character of $\mathbb{F}_q$, \[\mu(x) = \zeta^{\Tr (x)}\,.
\]

The Gauss sum (see \cite{LN,Wash}) of a character $\chi \in \widehat{\mathbb{F}_q^
{\times}}$  is defined as \[\tau(\chi) = -\sum_{x \in \mathbb{F}_q}\chi(x)\mu(x)\,.\]
We define \[g(j):=\tau(\omega^{-j})\,.\]

For any positive integer $j$, let $\wt_p(j)$ denote the $p$-weight of $j$, i.e.,
\[ \wt_p(j)=\sum_i j_i
\] 
where $\sum_i j_ip^i$ is the $p$-ary
expansion of $j$.

\subsection{Trace and similar objects}\label{traces}

Consider again the trace function $\Tr : \f{q} \to \f{p}$,
\[
\Tr (c) = c + c^p + c^{p^2} + \dots + c^{p^{n-1}}.
\]

We wish to generalise this definition to a larger class of finite field sums, which includes the usual trace function as a special case.
\bigskip

\begin{defn}\label{trace}
Let $p$ be a prime, let $n\ge 1$ be an integer and let $q=p^n$.
For any $S \subseteq \mathbb{Z}/(q-1)\mathbb{Z}$ satisfying $S^p=S$
where
$S^p := \{s^p \ | \ s \in S\}$,
we define the function $\tau_S:\f{q} \to \f{p}$ by
\[\tau_S(c) := \sum_{s\in S} c^s\,.\]
\end{defn}

\begin{defn}\label{lifttrace}
Let $p$ be a prime, let $n\ge 1$ be an integer and let $q=p^n$.
For any $S \subseteq \mathbb{Z}/(q-1)\mathbb{Z}$ satisfying $S^p=S$
where
$S^p := \{s^p \ | \ s \in S\}$,
we define the function $\tfr{\tau}_S:\f{q} \to \mathbb{Q}_p(\xi)$ by
\[\tfr{\tau}_S(c) := \sum_{s\in S} \omega^s(c)\,.\]
\end{defn}

\begin{rem}For the set $W = \{ p^i \ | \ i \in \{0,\dots,n-1\} \}$, $\tau_W$ is the usual trace function. \end{rem}

\begin{rem}By the definition of the Teichm\"{u}ller character, for any set $S$ 
we have $\tfr{\tau}_S \equiv \tau_S \pmod{p}$. 
Thus we may consider $\tfr{\tau}_S$ to be a \emph{lift} of $\tau_S$, and this explains the notation.
For the set $W$ defined in the previous remark, we let $\tfr{\Tr}$ denote the function $\tfr{\tau_W}$.
Sometimes we call $\tfr{\Tr}$ the lifted trace.
\end{rem}

Other than the set $W$, for the case $p = 3$, we will be particularly concerned with the following sets:
\begin{align*}
X&:=\{r \in \{0,\dots,q-2\}|r = 3^i+3^j\},\ (i, j \text{ not necessarily distinct})\\
Y&:=\{r \in \{0,\dots,q-2\}|r = 3^i+3^j+3^k, i,j,k \text{ distinct}\},\\
Z&:=\{r \in \{0,\dots,q-2\}|r = 2\cdot3^i+3^j, i\ne j\}.
\end{align*}

\subsection{Stickelberger's theorem and the Gross-Koblitz formula}\label{sect::Stickelberger}

Let $\pi$ be the unique ($p-1$)th root of $-p$ in $\mathbb{Q}_p(\xi,\zeta)$ satisfying \[\pi \equiv \zeta - 1 \pmod{\pi^{2}}\,.\]
Wan \cite{Wan}
noted that the following improved version of Stickelberger's theorem is a direct consequence of the Gross-Koblitz formula (Theorem \ref{GKformula}).
\bigskip

\begin{thm}\cite{Wan}\label{Stick+}
Let $1 \le j < q-1$ be an integer and let $j =j_0+j_1p+\dots + j_{n-1}p^{n-1}$. Then 
\[g(j) \equiv \frac{\pi^{\wt_p(j)}}{j_{0}!\cdots j_{n-1}!} \pmod{\pi^{\wt_p(j)+p-1}}\,.\]
\end{thm}

Stickelberger's theorem, as usually stated,  is the same congruence modulo $\pi^{\wt_p(j)+1}$.

We have (see \cite{GK, Rob}) that
$(\pi)$ is the unique prime ideal of $\mathbb{Q}_{p}(\zeta, \xi)$ lying above $p$. 
Since $\mathbb{Q}_{p}(\zeta, \xi)$ is an unramified extension of $\mathbb{Q}_{p}(\zeta)$,
which is a totally ramified (degree $p-1$) extension of $\mathbb{Q}_{p}$, it follows that $(\pi)^{p-1} = (p)$
and $\nu_p(\pi) = \frac{1}{p-1}$.
Here $\nu_p$ denotes the $p$-adic valuation.

Theorem \ref{Stick+} implies that
$\nu_{\pi}(g(j)) = \wt_p(j)$, and because 
$\nu_{p}(g(j)) =\nu_{\pi}(g(j)) \cdot \nu_p (\pi)$ we get 
\begin{equation}\label{val}
\nu_p(g(j)) = \frac{\wt_p(j)}{p-1}\,.\end{equation}

A generalisation of Stickelberger's theorem is the Gross-Koblitz formula.
\begin{thm}\label{GKformula}(Gross-Koblitz formula) \cite{GK}.

Let $1 \le j < q-1$ be an integer. Then
\[
g(j) = \pi^{\wt_p(j)}\prod_{i=0}^{n-1}\Gamma_{p}\left(\left\langle \frac{p^{i}j}{q-1} \right\rangle\right)
\]
where $\langle x \rangle$ is the fractional part of $x$, and $\Gamma_p$ is the p-adic gamma function.
\end{thm}

Our proof in Section \ref{MainProof} studies the $\pi$-adic expansion of the Kloosterman sum,
and uses the Gross-Koblitz formula to get information on the coefficients.

\subsection{The $p$-adic gamma function}
The $p$-adic gamma function $\Gamma_{p}$, introduced in \cite{Morita},
is defined over $\mathbb{N}$ by \[\Gamma_{p}(k) = (-1)^k\prod_{\substack{t<k\\(t,p) = 1}}t \,,\]
and extends to $\Gamma_{p} : \mathbb{Z}_{p} \to \mathbb{Z}_{p}$ according to Theorem \ref{GenWilson} below. 

The following are two classical results (they appear in \cite{disquisitiones}) which can be rephrased in terms of the $p$-adic gamma function. Theorem \ref{GenWilson} appears in this form in \cite{Morita}.

\begin{thm}[Wilson's theorem]\label{Wilson}
Let $p$ be an odd prime. Then 
\[
\Gamma_{p}(p-1) \equiv 1 \pmod{p}.
\]
\end{thm}

\begin{thm}[Generalised Wilson's theorem]\label{GenWilson}

Let $p$ be a prime, and suppose $x \equiv y \pmod{p^k}$ for some integer $k$. If $p^k \ne 4$, then
\[
\Gamma_{p}(x) \equiv \Gamma_p(y) \pmod{p^k}.
\]
\end{thm}

\subsection{Fourier coefficients}\label{sect::Fourier}

Recall that $\mu(x) = \zeta^{\Tr (x)}$.
The Fourier transform of a function $f:\mathbb{F}_q \to \mathbb{C}$ at $a \in \mathbb{F}_q$ 
is defined to be
\[\tfr{f}(a) = \sum_{x \in \mathbb{F}_q} f(x)\mu(ax)\,.\]
The complex number $\tfr{f}(a) $ is called the Fourier coefficient of $f$ at $a$.

Consider monomial functions defined by $f(x) = \mu(x^d)$. 
When $d=-1$ we have $\tfr{f}(a) = \mathcal{K}_{p^n}(a)$.
By
Fourier analysis \cite{KatzGaussSums, LL} we have for any $d$ 
\[
\tfr{f}(a) =\frac{q}{q-1}+ \frac{1}{q-1} \sum_{j=1}^{q-2}\tau(\bar{\omega}^j)\ \tau
(\omega^{jd})\ \bar{\omega}^{jd}(a) 
\]
and hence
\[
\tfr{f}(a) \equiv - \sum_{j=1}^{q-2}\tau(\bar{\omega}^j)\ \tau
(\omega^{jd})\ \bar{\omega}^{jd}(a) \pmod{q}\,.
\]

Putting $d=-1=p^n-2$, this congruence becomes

\begin{equation}\label{firstKL}
\K_q(a) \equiv - \sum_{j=1}^{q-2}(g(j))^2\ \omega^{j}(a) \pmod{q}.
\end{equation}

We will use this in Section \ref{GKmod27}.

\section{Proof of Theorem \ref{main}}\label{MainProof}

Moisio \cite{CertainValues} considered the reduction of the minimal polynomial $m_a(x)$ modulo $p$,
and proved the following.

\begin{lem}\label{MinPolModp}\cite{CertainValues}
For $a \in \f{q}$, let $m(x)$ be the minimal polynomial of $\K_q(a)$ over $\QQ$ and let $t$ be the degree of $m$. Then  \[m(x) \equiv x^t \pmod{p}.\]
\end{lem}

Our first result concerns the reduction of the minimal polynomial $m_a(x)$ modulo $p^2$.
  
\bigskip

\begin{thma}\label{main2}
Let $p$ be an odd prime, and let $\left(\frac{\cdot}{p}\right)$ be the Legendre symbol. Then 
\[\prod_{i=1}^{\frac{p-1}{2}}(\K_q(i^2a)) \equiv p \left(\frac{\Tr(a)}{p}\right) \pmod{p^2}\,.\]
\end{thma}

Proof: For $j \in \{1,\dots,q-2\}$, Theorem \ref{Stick+} implies that 
\begin{equation}\label{pi_val}
\nu_{\pi}(g(j)^2) = 2\wt_p(j), 
\end{equation}
so taking equation \eqref{firstKL} mod ${\pi^4}$ gives
\begin{align*}
\K_q(a) &\equiv - \sum_{\wt_p(j) = 1}g(j)^2\ \omega^{j}(a) \pmod{\pi^4}\\
&\equiv -g(1)^2\tfr{\Tr}(a) \pmod{\pi^4}.
\end{align*}

Equation \eqref{pi_val} implies that $\nu_{\pi}(g(1)^2) = 2$. Therefore we can write $\K_q(a)$ as 
\[
\K_q(a) = a_1\pi^2+a_2\pi^4+\cdots ,
\]
where 
\begin{align*}
a_1 &= -\left(\frac{g(1)}{\pi}\right)^2\tfr{\Tr}(a)\\
& = -\left(\prod_{i=0}^{n-1}\Gamma_{p}\left(\left\langle \frac{p^{i}}{q-1} \right\rangle\right)\right)^2\tfr{\Tr}(a) \text{ (by Theorem \ref{GKformula}).}
\end{align*}

Reducing this expression modulo $p$ gives that

\begin{align*}
a_1 &\equiv -\left(\Gamma_{p}\left(\frac{1}{q-1}\right)\right)^{2}\Tr(a) \pmod{p}\\
&\equiv -\left(\Gamma_p(p-1) \right)^{2}\Tr(a) \pmod{p} \text{ (by Theorem \ref{GenWilson})}\\
&\equiv-\Tr(a) \pmod{p} \text{ (by Theorem \ref{Wilson}),}
\end{align*}
and thus
\[
\K_q(a) \equiv -\pi^2\Tr(a) \pmod{\pi^4}.
\]

So
\begin{align*}
\prod_{i=1}^{\frac{p-1}{2}}(\K_q(i^2a)) \equiv& \pi^{p-1}\prod_{i=1}^{\frac{p-1}{2}}(-i^2\Tr(a)) \pmod{\pi^{p+1}}\\
&\\
\equiv &-p \Tr(a)^{\frac{p-1}{2}}\prod_{i=1}^{\frac{p-1}{2}}(-i^2) \pmod{\pi^{p+1}}.
\end{align*}
But $\prod_{i=1}^{\frac{p-1}{2}}(\K_q(i^2a)) \in \ZZ$ by the remarks in Section \ref{Intro}, so
\[
\prod_{i=1}^{\frac{p-1}{2}}(\K_q(i^2a)) \equiv -p \Tr(a)^{\frac{p-1}{2}}\prod_{i=1}^{\frac{p-1}{2}}(-i^2) \pmod{p^2}.
\]
Using Wilson's Theorem (as usually stated), we have that
\[
\prod_{i=1}^{\frac{p-1}{2}}(-i^2) =\prod_{i=1}^{p-1} i \equiv -1 \pmod{p}.
\]
Thus
\[
\prod_{i=1}^{\frac{p-1}{2}}(\K_q(i^2a)) \equiv p\Tr(a)^{\frac{p-1}{2}} = p\left(\frac{\Tr(a)}{p}\right)\pmod{p^2}.
\]\qed

\bigskip

\begin{cor}
The constant term of the characteristic polynomial $c_a(x)$
is always congruent to either 0 or $\pm p$ mod $p^2$.
\end{cor}

The following result is due to Wan.
\begin{thm}\cite{Wan}
Let $a \in \f{q}$. If $\Tr(a) \ne 0$, the minimal polynomial of $\K_q(a)$ has degree $\frac{p-1}{2}$.
\end{thm}

Thus if $\Tr(a) \ne 0$, the minimal polynomial $m(x)$ of $\K_q(a)$ is precisely the characteristic polynomial $c(x)$. In this case (and in the case that $\deg(m(x)) = \frac{p-1}{2}$ where $\Tr(a) = 0$) Theorem \ref{main} gives a statement about the constant term of $m(x)$ mod $p^2$.

If $\Tr(a) = 0$ and $\deg(m(x)) <\frac{p-1}{2}$, then the result in Theorem \ref{main} is implied by Lemma \ref{MinPolModp}. In this case, our result gives us no extra information about the constant term of the minimal polynomial.

\section{Ternary Kloosterman sums modulo 27}\label{GKmod27}

In this section we use the same techniques to improve the modulo 9 Kloosterman sum 
characterisation in \cite{SetaFGR}
to a modulo 27 characterisation.
First let us prove a lemma on evaluations of the $p$-adic gamma function. This lemma will allow us to 
evaluate Gauss sums for higher moduli and find Kloosterman congruences modulo $27$.
\bigskip

\begin{lem}\label{gammamod}Let $n \ge 3$ $q=3^n$ and let $i$ be an integer in the range $0,\dots n-1$. Then
\[
\Gamma_3\left(\left\langle \frac{3^{i}}{q-1} \right\rangle\right) \equiv 
\left\{ \begin{array}{ll}
13 \pmod{27} & \textrm{if } i = 1, \\
1 \pmod{27} & \textrm{if } i > 1.
\end{array} \right.
\]
\end{lem}

\begin{proof}
For any $3 \le j \le n$, we have $3^j \le q$, and
\[
\left\langle\frac{3^i}{q-1}\right\rangle = \frac{3^i}{q-1} \equiv 3^i(3^j-1) \pmod{3^j},\] so 
\[
\Gamma_3\left(\left\langle \frac{3^{i}}{q-1} \right\rangle\right) \equiv \Gamma_3(26\cdot3^i) \pmod{27}.\]
If $i \ge 3$, then $26\cdot 3^i \equiv 0 \pmod{27}$, and
\[ 
\Gamma_3\left(\left\langle \frac{3^{i}}{q-1} \right\rangle\right) \equiv 1 \pmod{27}\,,
\]

Now $\Gamma_3(26\cdot3) \equiv \Gamma_3(24) \pmod{27}$ using Theorem \ref{GenWilson}. And $\Gamma_3(24) \equiv 13 \pmod{9}$. Similarly:
\begin{eqnarray*}
\Gamma_3(26\cdot9) & \equiv & 1 \pmod{27}.
\end{eqnarray*}
\end{proof}

Lemma \ref{gammamod} allows us to compute Gauss sums modulo $27$:
\bigskip

\begin{lem}\label{wt1lem}
Let $n \ge 3$ and let $q=3^n$. Then
\[
g(j)^2 \equiv \left\{ \begin{array}{ll}
6 \pmod{27} & \textrm{if } \wt_p(j) =  1, \\
9 \pmod{27} & \textrm{if } \wt_p(j) =  2, \\
0 \pmod{27} & \textrm{if } \wt_p(j) \ge 3.
\end{array} \right.
\]
\end{lem}

\begin{proof}
Suppose $\wt_3(j) = 1$. By Theorem \ref{GKformula} and Lemma \ref{gammamod}, 
\[
g(j) \equiv 13 \pi \pmod{27}.
\]
Let 
\[
g(j) = 27A + 13\pi
\]
for some $A\in \mathbb{Z}_{3}[\zeta, \xi]$. Then
\begin{align*}
g(j)^2 & = 27^2 A^2 + 2\cdot 27\cdot 13 A  + 169\pi^{2} \\
			 & \equiv 169\pi^{2} \pmod{27} \\
			 & \equiv 6 \pmod{27}
\end{align*}
since $\pi^2=-3$. Now suppose $\wt_3(j) = 2$. By Theorem \ref{GKformula}, 
\[
g(j) \equiv -3 \pmod{9}. 
\] Thus $g(j) = 9B-3$ for some $B\in \mathbb{Z}_{3}[\zeta, \xi]$,
so \[g(j)^2 = 81B^2-54B+9 \equiv 9 \pmod{27}.\]

It is clear from Theorem \ref{GKformula} that if $\wt_3(j)>2$, then
\[
27 | \pi^{2\wt_3(j)} | g(j)^2.
\]
\end{proof}

Now we are ready to prove our result on Kloosterman sums modulo $27$.
\bigskip

\begin{thm}\label{mod27}Let $n \ge 3$, $q=3^n$ and
let $\tfr{\Tr}$ and $\tfr{\tau_{X}}$ be as defined in Section \ref{traces}.
Then

\begin{equation}\label{hateq}
\K_{3^n}(a) \equiv 21\tfr{\Tr}(a)+18\tfr{\tau_{X}}(a) \pmod{27}.
\end{equation}
\end{thm}
\begin{proof}
Using   \eqref{firstKL} and Lemma \ref{wt1lem}, we get
\begin{align*}
\mathcal{K}(a) & \equiv - \sum_{j=1}^{q-2}g(j)^2\ \omega^{j}(a)\pmod{q} \\
		& \equiv - \sum_{\wt_{3}(j)=1}g(j)^2 {\omega}^{j}(a) - \sum_{\wt_{3}(j)=2} g(j)^2 {\omega}^{j}(a) \pmod{27} \\
		& \equiv -6\sum_{\wt_{3}(j)=1}{\omega}^{j}(a) - 9\sum_{\wt_{3}(j)=2}{\omega}^{j}(a) \pmod{27} \\
&\equiv 21\tfr{\Tr}(a)+18\tfr{\tau_{X}}(a) \pmod{27}.\qedhere
 \end{align*}
\end{proof}

Next we shall  express the above result in terms of operations within $\f{q}$ itself, i.e.,
using functions $\tau_S$ directly, and not their lifts.
Note that in (\ref{hateq}) we only need $\tfr{\Tr}(a)$ modulo $9$ and 
$\tfr{\tau_{X}}(a)$ modulo $3$. We have
\[
\tau_{X} (a) \equiv \tfr{\tau_{X}}(a) \pmod{3}
\]
so this takes care of the $\tfr{\tau_{X}}(a)$ term.
For the other term
we need to find a condition for $\tfr{\Tr}(a)$ modulo $9$ using functions from $\f{q}$ to $\f{3}$. We will do that in the proof of the following corollary.
\bigskip

\begin{cor}\label{fieldsums}
Let $n \ge 3$, $q=3^n$, $a\in\f{q}$ and
let $\tau_{X}$, $\tau_{Y}$ and $\tau_{Z}$ be as defined in Section \ref{traces}. Let $\Tr(a)$ be the trace of $a$, but considered as an integer. Then 
\[\K_q (a) \equiv 21 \Tr(a)^3+18\tau_{Z}(a)+9\tau_{Y}(a)+18\tau_{X}(a) \pmod{27}.\]
\end{cor}
\begin{proof}

First recall that $\tfr{\tau}_{X}(a) \equiv \tau_{X}(a) \pmod{3}$.

To determine $\tfr{\Tr}(a) \bmod{9}$, we compute 
\begin{align*}
\tfr{\Tr}(a)^3 &= \sum_{i,j,k \in \{0,\dots,n-1\}}\omega(a^{3^i+3^j+3^k})\\
&=\tfr{\Tr}(a)+3\tfr{\tau}_Z(a)+6\tfr{\tau}_Y(a)\,,
\end{align*}
and note the elementary fact that if $x \equiv y \pmod{m}$, then $x^m \equiv y^m \pmod{m^2}$. This means that $\tfr{\Tr}(a)^3 \bmod{9}$ is given by $\tfr{\Tr}(a) \bmod{3} = \Tr (a)$, i.e. $\tfr{\Tr}(a)^3 \bmod{9} = \Tr(a)^3$.

Since \[\tfr{\tau}_Z(a) \equiv \tau_{Z}(a)\pmod{3}\]and \[\tfr{\tau}_Y(a) \equiv \tau_{Y}(a)\pmod{3}\,,\]we have that

\[\tfr{\Tr}(a) \equiv \Tr(a)^3 - 3\tau_{Z}(a) - 6\tau_{Y}(a)\pmod{9},\]proving the result.
\end{proof}

The next corollary combines Corollary \ref{fieldsums} and Theorem \ref{mod27} and
enumerates the possible values of ternary Kloosterman sums mod 27. 

\begin{cor}\label{main27}
Let $n\ge 3$, and let $q=3^n$. Let $\Tr$, $\tau_{X}$ and $\tau_{Y}$ be as defined in Section \ref{traces}.  Then 
\begin{equation*}
\K_q(a)\equiv \left\{
	\begin{array}{rcccccl}
		0\pmod{27}\text{ if } &\Tr(a) = &0 &\text{ and } &\tau_{Y}(a) &+2\tau_{X}(a)&=0\\
		3\pmod{27}\text{ if } &\Tr(a) = &1 &\text{ and } &\tau_{Y}(a) & &= 2\\
		6\pmod{27}\text{ if } &\Tr(a) = &2 &\text{ and } &\tau_{Y}(a)&+\tau_{X}(a) &= 2\\
		9\pmod{27}\text{ if } &\Tr(a) = &0 &\text{ and } & \tau_{Y}(a)&+2\tau_{X}(a) &= 1\\
		12\pmod{27}\text{ if } &\Tr(a) = &1 &\text{ and } & \tau_{Y}(a)& &= 0\\
		15\pmod{27}\text{ if } &\Tr(a) = &2 &\text{ and } &\tau_{Y}(a)&+\tau_{X}(a) &= 0\\
		18\pmod{27}\text{ if } &\Tr(a) = &0 &\text{ and } &\tau_{Y}(a)&+2\tau_{X}(a) &=2\\
		21\pmod{27}\text{ if } &\Tr(a) = &1 &\text{ and } &\tau_{Y}(a)& &= 1\\
		24\pmod{27}\text{ if } &\Tr(a) = &2 &\text{ and } &\tau_{Y}(a)&+\tau_{X}(a) &= 1.
		\end{array} \right.
\end{equation*}
\end{cor}

\begin{proof}
Note that
\[\Tr(a) \tau_{X}(a) = \Tr(a)+2\tau_{Z}(a)\,.\]
Thus Corollary \ref{fieldsums} can be rewritten as
\begin{equation}\label{final27}
\K_q (a) \equiv 21 \Tr(a)^3+18\Tr(a)+18\tau_{X}(a)+9\Tr(a)\tau_{X}(a)+9\tau_{Y}(a) \pmod{27}.
\end{equation}

The result is an enumeration of the cases in equation \eqref{final27}.\qedhere
\end{proof}

\bigskip

We remark that a characterisation like in Corollary \ref{main27} of Kloosterman sums
modulo $p^3$ for $p>3$ does not seem to be straightforward.  
The estimates given by the Gross-Koblitz formula are weaker.


\bibliographystyle{plain}
\bibliography{Kloosterman}
\end{document}